\documentclass[11pt]{amsart}
\usepackage{graphicx, color, amssymb}

\newtheorem{Theorem}{Theorem}

\newtheorem{Definition}[Theorem]{Definition}

\begin{document}

\title[Complex of non-crossing diagonals]{The complex of non-crossing diagonals of a polygon}

\author[B. Braun]{Benjamin Braun}

\address{Department of Mathematics, University of Kentucky, Lexington, KY}
\urladdr{http://www.ms.uky.edu/$\sim$braun}
\email{braun@ms.uky.edu}

\author[R. Ehrenborg]{Richard Ehrenborg}
\thanks{The second author was partially supported by National Security Agency grant H98230-06-1-0072.}
\address{Department of Mathematics, University of Kentucky, Lexington, KY}
\urladdr{http://www.ms.uky.edu/$\sim$jrge}
\email{jrge@ms.uky.edu}

\date{February 10, 2008}

\begin{abstract}
Given a convex $n$-gon $P$ in the Euclidean plane, it is well known that the simplicial complex $\theta(P)$ with vertex set given by diagonals in $P$ and facets given by triangulations of $P$ is the boundary complex of a polytope of dimension $n-3$.  We prove that for any non-convex polygonal region $P$ with $n$ vertices and $h+1$ boundary components, $\theta(P)$ is a ball of dimension $n+3h-4$.  We also provide a new proof that $\theta(P)$ is a sphere when $P$ is convex.
\end{abstract}

\maketitle

An $n$-gon (or \textit{polygon}) $P$ in the Euclidean plane is the figure formed by $n$ distinct points $x_1,\ldots,x_n$ in the plane, called \textit{vertices}, and $n$ line segments $[x_i,x_{i+1}]$, $i=1,\ldots,n$, called \textit{edges}, with addition modulo $n$ on the indices, subject to the following condition: the only points of the plane that belong to two edges of $P$ are the vertices of $P$.  It is well know that the edges of $P$ form a Jordan curve in the plane and thus $P$ has a well defined interior and exterior.  We view the polygon $P$ as consisting of both its interior and its boundary.  A \textit{diagonal} in $P$ is a line segment $[x,y]$ between two non-adjacent vertices of $P$ that is contained in $P$.

One may also consider more general \textit{polygonal regions} in the plane, i.e. bounded, connected regions whose boundary is the disjoint union of the boundaries of polygons.  Vertices, edges, and diagonals are defined analogously to that for polygons.  It is well known that for any polygonal region with $h+1$ boundary components, one can triangulate $P$ using $n+3h-3$ diagonals of $P$.  Triangulations of polygons are of great interest in discrete and computational geometry, for example see \cite{AlgorGeom}.  We are interested in the following topological structure on the set of diagonals in $P$.

\begin{Definition} For any polygonal region $P$, let {\em the complex of non-crossing diagonals of $P$}, $\theta(P)$, be the simplicial complex with vertex set the diagonals in~$P$ and facets given by triangulations of $P$.
\end{Definition}

Since any triangulation of $P$ uses $n+3h-3$ diagonals, we see that $\theta(P)$ is a pure simplicial complex, i.e. all the maximal faces of $P$ have the same dimension.  It was further shown in the late 1980's, independently by M.\ Haiman \cite{HaimanAssociahedron} and C.\ Lee \cite{LeeAssociahedron}, that for convex polygons $\theta(P)$ arises as the boundary complex of a polytope of dimension $n-3$ called the \textit{associahedron}.  Thus, $\theta(P)$ triangulates a sphere of dimension $n-4$.  Associahedra arise from various constructions and have been studied intensely, see \cite{ZieglerLectures} and the references therein for further discussion.

Our purpose in this note is to prove the following theorems on the structure of $\theta(P)$ for polygons and polygonal regions.

\begin{Theorem}\label{nonconvex}
For any non-convex $n$-gon, $\theta(P)$ is homeomorphic to $B^{n-4}$, a ball of dimension $n-4$.  For any convex $n$-gon $P$, $\theta(P)$ is homeomorphic to $S^{n-4}$, a sphere of dimension $n-4$.
\end{Theorem}

\begin{Theorem}\label{region}For any polygonal region $P$ with $n$ vertices and $h+1$ boundary components, $\theta(P)$ is homeomorphic to $B^{n+3h-4}$, a ball of dimension $n+3h-4$.  Further, every vertex of $\theta(P)$ is on the boundary of $\theta(P)$.
\end{Theorem}

\begin{figure}[ht]
\begin{center}
 \input{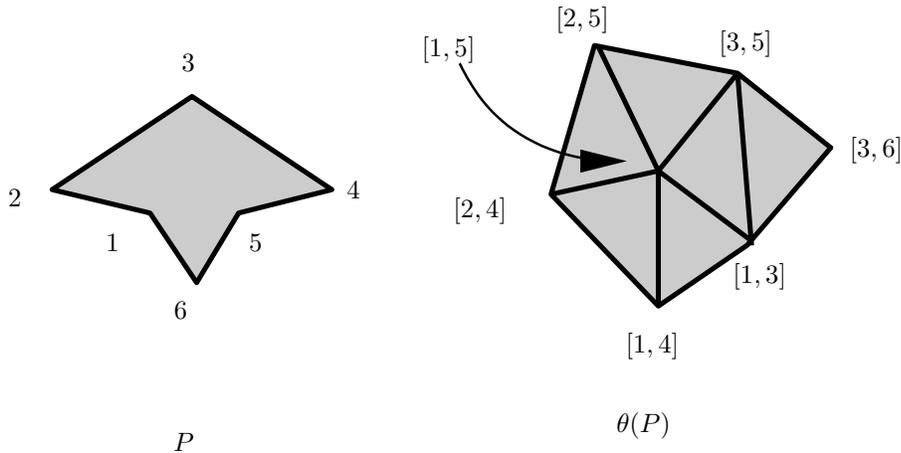}
\end{center}
\caption{A non-convex polygon $P$ and its complex of non-crossing diagonals $\theta(P)$.}
\end{figure}

Our proof that $\theta(P)$ is spherical for convex $P$ is new, though not as strong as the polytopal result mentioned earlier.  We present it here because we have found it illuminating regarding the spherical structure of $\theta(P)$ and hope it will be so for others.

Motivating examples for our general method can be found by looking at small convex polygons.  Let $P_n$ represent a convex $n$-gon with vertices labeled $1$ through $n$.  It is clear that $\theta(P_4)$ consists of two points, i.e. a zero dimensional sphere.  We obtain $P_5$ by inserting $5$ between $1$ and $4$, introducing three new diagonals: $[1,4],[2,5]$ and $[3,5]$.  These diagonals form two disjoint faces in $\theta(P_5)$ consisting of the first diagonal and the last two.  The diagonal $[1,4]$ may be added to any triangulation of $P_4$, producing a cone over $\theta(P_4)$ in $\theta(P_5)$.  The $1$-cell  $\{[2,5],[3,5]\}$ in $\theta(P_5)$ is connected on its boundary points to $[1,3]$ and $[2,4]$, respectively.  It is clear that the midpoint of $\{[2,5],[3,5]\}$ acts as another cone point over $\theta(P_4)$ in $\theta(P_5)$, thus showing that $\theta(P_5)$ is the suspension of $\theta(P_4)$.  The case $n$ equal to $6$ can be analyzed similarly, with two cone points over $\theta(P_5)$ coming from the diagonal $[1,5]$ and the barycenter of the $2$-cell $\{[2,6],[3,6],[4,6]\}$.

In general, when we pass from an $n$-gon to an $(n+1)$-gon, we are adding in $n-1$ new diagonals which are naturally partitioned into two sets: the $n-2$ diagonals adjacent to the new vertex $n+1$ and the single diagonal $[1,n]$.  Topologically, adding these two sets of diagonals to $\theta(P_n)$ yields a suspension, producing a spherical structure.  The vertex $n+1$ plays a special role for this analysis, which is captured by the second half of the following definition.

\begin{Definition}{\rm (\cite{Meisters}, \cite{Toussaint})} A vertex $x_i$ of a polygon $P$ is called a {\em principal vertex} provided the associated line segment $[x_{i-1},x_{i+1}]$ intersects the boundary of $P$ only at $x_{i-1}$ and $x_{i+1}$.  A principal vertex $x_i$ of $P$ is called a {\em mouth} if the associated line segment $[x_{i-1},x_{i+1}]$ is external to $P$, i.e., the interior of $[x_{i-1},x_{i+1}]$ lies in the exterior of $P$.  A principal vertex $x_i$ of $P$ is called an {\em ear} if the associated line segment $[x_{i-1},x_{i+1}]$ is interior to $P$, i.e., the interior of $[x_{i-1},x_{i+1}]$ lies in the interior of $P$.
\end{Definition}

\begin{Theorem}\label{ears}{\rm (G. Meisters, \cite{Meisters})} Every polygon has at least two ears.
\end{Theorem}

\begin{Theorem}\label{mouth}{\rm (G. Toussaint, \cite{Toussaint})} Every non-convex polygon has at least one mouth.
\end{Theorem}

An important point is that convex polygons have only ears and no mouths.  In the case when $P$ is not a convex polygon, then $P$ has at least one mouth; these special vertices end up forcing $\theta(P)$ to be collapsible.  To demonstrate the collapsing, we will use the discrete Morse theory developed by R. Forman in \cite{FormanMorseTheory}.  Specifically, we will use a simplified version of the ``Pairing Lemma'' found in the work \cite{ShareshianLinusson} of J.\ Shareshian and S.\ Linusson.

\begin{Theorem}\label{Pairing}{\rm (J.\ Shareshian and S.\ Linusson, \cite{ShareshianLinusson})} Let $\Sigma$ be a simplicial complex on a partially ordered vertex set $(V,\preceq)$.  Let $Q$ denote the face poset of $\Sigma$, including $\emptyset$.  For a function $f:Q\rightarrow V$, set $Q_f:=\{\sigma\in Q: f(\sigma)\notin \sigma\}.$  For $\sigma \in Q_f$, set $\sigma^+:=\sigma\cup\{f(\sigma)\},$ and for $\tau \in Q \backslash Q_f$, set $\tau^-:=\tau\backslash\{f(\tau)\}.$  Assume that $f$ satisfies the following conditions:
\begin{enumerate}
\item{If $\sigma\in Q_f$, then $\sigma^+\in Q$.}
\item{If $\sigma \in Q_f$, then $f(\sigma^+)=f(\sigma)$.}
\item{If $\tau \in Q\backslash Q_f$ and $\tau^-\in Q$, then $f(\tau^-)=f(\tau)$.}
\item{If $x\in \sigma\in Q_f$ and $\sigma^+\backslash \{x\}\in Q$, then $f(\sigma)\preceq f(\sigma^+\backslash \{x\})$.}
\end{enumerate}
Then the simplicial complex $\Sigma$ is collapsible.
\end{Theorem}

For the situation described above, the original statement of Linusson and Shareshian provides only contractibility of $\Sigma$, but Forman's general theory implies the stronger collapsibility condition; details may be found in \cite{FormanMorseTheory}.  To show that $\theta(P)$ is actually a ball in the non-convex cases of our proof below, we will apply the following theorem originally due to Whitehead.

\begin{Theorem}\label{ball}{\rm (Theorem 1.6 of \cite{FormanMorseTheory})} Let $M$ be a piecewise linear $n$-manifold with boundary and $x$ a vertex of $M$.  If $M$ collapses to $x$, then $M$ is a piecewise linear $n$-ball.
\end{Theorem}

\textbf{Proof of Theorem \ref{nonconvex}:} We proceed by induction on $n$.  Assume $n\geq 5$, as the base case $n=4$ is clear.  For a vertex $x$ of $P$, let \[D_{x}:=\{[x,w_1],\ldots,[x,w_k]\}\] be the set of diagonals in $P$ adjacent to $x$.  Suppose that $P$ is a non-convex $n$-gon and let $x$ be a mouth of $P$.  Given a triangulation $T$ of $P$, $x$ is contained in some triangle $\Delta$ in $T$.  If no diagonal adjacent to $x$ is an edge of $\Delta$, then $x$ must be an ear, contradicting our assumption.  Thus, for any triangulation $T$ of $P$, $T$ contains a diagonal incident to $x$.

Linearly order the diagonals in $P$ in such a way that the last $k$ elements in the ordering are $[x,w_1],\ldots,[x,w_k]$. We will construct a function $f$ from the face poset of $\theta(P)$ to the (ordered) set of diagonals in $P$. Any face $\sigma$ of $\theta(P)$ consists of a set of pairwise non-crossing diagonals in $P$.  Any such set can be extended to a triangulation of $P$, hence there is some diagonal in $D_{x}$ which is pairwise non-crossing with the elements of $\sigma$.  Set $f(\sigma)$ equal to the diagonal $[x,w_j]$ with the property that $j$ is maximized over all elements of $D_{x}$ that are pairwise non-crossing with the elements of $\sigma\backslash D_{x}$.  It is an easy verification that our $f$ satisfies conditions $(1)$ through $(4)$ in Theorem \ref{Pairing}.

To conclude that $\theta(P)$ is a ball, we only need to check that $\theta(P)$ is a piecewise linear manifold.  Given any vertex $v$ in $\theta(P)$, i.e. a diagonal in $P$, we check that the link of $v$ is either a sphere or a ball.  It is clear that $v$ cuts $P$ into two polygons $G$ and $H$ with $j$ and $n-j+2$ vertices, respectively, for some $j$.  If $G$ and $H$ are both convex, then by induction $\theta(G)$ and $\theta(H)$ are spheres of dimension $j-4$ and $n-j-2$, respectively.  The link of $v$ in $\theta(P)$ in this case is $\theta(G)\ast\theta(H)\cong S^{j-4} \ast S^{n-j-2}\cong S^{n-5}$, where $\ast$ denotes the join operation, see \cite{BjornerSurvey} and \cite{Hatcher} for definitions.  Without loss of generality, if $G$ is convex and $H$ is non-convex, then by induction the link of $v$ in $\theta(P)$ is $\theta(G)\ast\theta(H)=S^{j-4}\ast B^{n-j-2}\cong B^{n-5}$.  Finally, if both $G$ and $H$ are non-convex, we see that by induction the link of $v$ in $\theta(P)$ is $\theta(G)\ast\theta(H)=B^{j-4}\ast B^{n-j-2}\cong B^{n-5}$.  Thus, $\theta(P)$ is a piecewise linear manifold and Theorem \ref{ball} applies.

For the case where $P=P_n$ is a convex $n$-gon, consider $\theta(P_n)$ as the union of two closed subspaces $\mathcal{C}:=\theta(P_{n-1})\ast[1,n-1]$ and $\mathcal{J}:=\mathrm{dl}_{\theta(P_{n})}([1,n-1])$, where $\mathrm{dl}_{\Sigma}(\sigma)$ is the deletion of $\sigma$ from $\Sigma$, see again \cite{BjornerSurvey} for definitions.  Note that $\mathcal{C}\cap\mathcal{J}=\theta(P_{n-1})$ and, by our induction hypothesis, this is a sphere of dimension $n-5$.  The complex $\mathcal{C}$ is a cone over $\theta(P_{n-1})$ and hence is a ball of dimension $n-4$.  The complex $\mathcal{J}$ is isomorphic to the complex $\theta(R)$ obtained by adding an extra vertex to $P_{n-1}$ just inside the edge $[1,n-1]$ to obtain a non-convex polygon $R$.  By induction, this is also a ball of dimension $n-4$ with boundary easily observed to be $\theta(P_{n-1})$.  Hence, $\theta(P_n)$ is the union of two $(n-4)$-balls glued homeomorphically along their boundary spheres and is itself a sphere of dimension $n-4$.

\textbf{Proof of Theorem \ref{region}:} Let $P$ be a polygonal region with $n$ vertices and $h+1$ boundary components.  Denote the boundary of $P$ as the union $C_1\cup\cdots\cup C_{h+1}$ of the boundary components $C_i$ of $P$, where $P\subseteq \mathrm{conv}(C_1)$, i.e. $C_1$ is the outermost bounding curve for $P$.  Consider a mouth $x$ of $P$, obtained as a mouth of $C_1$ or an ear of $C_i$ for some $i\neq 1$.  An identical analysis to the case where $h=0$ shows that any triangulation of $P$ contains a diagonal adjacent to $x$ and that the function assigning to each face $\sigma$ of $\theta(P)$ the maximal $j$ in $D_x$ such that $[x,w_j]$ may be added to $\sigma\backslash D_x$ yields via Theorem \ref{Pairing} a collapsing to a point.  What remains is only to check that $\theta(P)$ is  a combinatorial manifold, enabling us to invoke Theorem \ref{ball}.

This situation is different from our earlier case in that a diagonal in $P$ need not disconnect the interior of $P$.  If it does, then by induction the link of the diagonal is the join of either a ball and a sphere or of two balls, and hence is a ball as desired.  If not, let $[x,w]$ denote the diagonal in question.  We may add two new vertices $x'$ and $w'$ overlapping $x$ and $w$, respectively, and a new diagonal $[x',w']$ overlapping $[x,w]$.  By slightly perturbing $x$ and $w$ to lie outside of $P$, we obtain a new polygon $P'$ from $P$ with two new vertices and two new edges, but with one less boundary component, namely we merge the components $C_w$ and $C_x$ containing $w$ and $x$, respectively, along the diagonals $[x,w]$ and $[x',w']$.  We may perturb these vertices in such a way that $P'$ has the same set of diagonals as those in $P$ not crossing $[x,w]$, i.e. $\theta(P')$ is homeomorphic to the link of $[x,w]$ in $\theta(P)$.  By induction, we conclude that $\theta(P')$ is a ball of dimension $(n+2)+3(h-1)-4=n+3h-5.$ As the link of every vertex in $\theta(P)$ is a ball, Theorem \ref{ball} may be applied.  Further, this demonstrates that every vertex of $\theta(P)$ lies on the boundary of $\theta(P)$.

\bibliographystyle{plain}
\bibliography{Braun}

\end{document}